\documentclass[a4paper,10pt,leqno]{article}
\usepackage{epsfig, graphics, amssymb, amsfonts, euscript ,stmaryrd,float,graphicx,color}
\makeatletter
       \def\ps@plain{
      \def\@oddhead{}\def\@oddfoot{\rm\hfil LAMFA,  No.9/Mars 2004, p. \thepage}
       \def\@evenhead{}\let\@evenfoot\@oddfoot}
       \makeatother
\def\today{\number\day\space\ifcase\month\or
January\or February\or March\or April\or May\or June\or July\or
August\or September\or October\or November\or December\fi
\space\number\year}

\makeatletter

\def\@begintheorem#1#2{\trivlist \item[\hskip \labelsep{\sc #1\ #2.}]\sl}
\def\@opargbegintheorem#1#2#3{\trivlist
\item[\hskip \labelsep{\sc #1\ #2\ (#3).}]\sl}
\def\@endtheorem{\endtrivlist}

\def\@sect#1#2#3#4#5#6[#7]#8{\ifnum #2>\c@secnumdepth
\let\@svsec\@empty\else
\refstepcounter{#1}\edef\@svsec{\csname the#1\endcsname.\hskip
0.5em}\fi \@tempskipa #5\relax \ifdim \@tempskipa>\z@
\begingroup #6\relax
\@hangfrom{\hskip #3\relax\@svsec}{\interlinepenalty \@M #8\par}
\endgroup \csname #1mark\endcsname{#7}\addcontentsline
{toc}{#1}{\ifnum #2>\c@secnumdepth \else \protect\numberline{\csname
the#1\endcsname}\fi #7}\else
\def\@svsechd{#6\hskip #3\relax\@svsec #8.\csname #1mark\endcsname
{#7}\addcontentsline {toc}{#1}{\ifnum #2>\c@secnumdepth \else
\protect\numberline{\csname the#1\endcsname}\fi #7}}\fi \@xsect{#5}}

\makeatother 
\newtheorem{Th}{Theorem}[section]
\newtheorem{Prop}{Proposition}[section]
\newtheorem{Lm}{Lemma}[section]

\newtheorem{Rm}{Remark}[section]

\newenvironment{Proof}{\removelastskip\vskip12pt plus 1pt \noindent\sc
Proof.\quad\rm}{\hspace*{\fill}$\blacksquare$}
\newenvironment{Pf}{\removelastskip\vskip12pt plus 1pt \noindent
\rm}{\hspace*{\fill}$\blacksquare$}

\newcommand{\be}{\begin{equation}}
\newcommand{\bT}{\begin{array}}
\newcommand{\eT}{\end{array}}
\newcommand{\de}{\end{equation}}
\newcommand{\ee}{\end{equation}}

\def\t{\tau}

\def\text#1{{\quad \hbox{#1} \quad}}

\def\ds{\displaystyle}
\def\sec{{\prime\prime}}

\begin{document}
\title{\bf On  similarity and pseudo-similarity solutions of  Falkner-Skan boundary layers}
\author{
 \bf M.~\textsc{\bf Guedda} and  \bf Z.~\textsc{\bf Hammouch}\footnote{Corresponding author. E-mail: zakia.hammouch@u-picardie.fr}  \\
{\footnotesize \it LAMFA, CNRS UMR 6140, Universit\'e de Picardie
Jules Verne,
}\\
{\footnotesize \it Facult\'e de Math\'ematiques et d'Informatique,
33, rue Saint-Leu 80039 Amiens, France} }
\date{}
\maketitle

\noindent{\bf Abstract.} This paper deals with the two-dimensional
incompressible, laminar, steady-state boundary layer
equations.First, we determine a family of velocity distributions
outside the boundary layer such that these problems may have
similarity solutions.Then, We examine in detail new exact solutions,
called {\it Pseudo-similarity}, where the external velocity varies
inversely-linear with the distance in the $x$-direction along the
surface ($U_e(x) = U_\infty x^{-1}). $ The analysis shows that
solutions exist only for a lateral suction. Here it is assumed that
the flow is induced by a continuous permeable surface with the
stretching velocity  $ U_wx^{-1}.$ For specified conditions, we
establish the existence of  an infinite number of solutions,
including monotonic solutions and solutions which oscillate an
infinite number of times and tend to a certain
 limit. The properties of solutions depend on the
suction parameter. Furthermore, making use of the fourth-order
Runge-Kutta scheme together with the shooting method, numerical
solutions are obtained.\\

\noindent {\it keywords}: Boundary layer,
Falkner-Skan, Similarity solution, Pseudo-similarity.


\section{Introduction}\label{Sect.1}
\setcounter{equation}{0} \setcounter{theorem}{0}
\setcounter{lemma}{0} \setcounter{remark}{0}
\setcounter{corollary}{0}

\qquad In this paper we are concerned with the classical
two-dimensional laminar incompressible boundary layer flow past a
wedge or a flat plate \cite{Sch}. For the first approximation, the
model  is described by   the Prandlt equations or the boundary layer
equations

\begin{equation}\label{eq:Prandtl}
u\partial_x u + v\partial_y u= U_e\partial_x U_e+\nu\partial^2_{yy}
u, \quad \partial_x u + \partial_y v = 0,
\end{equation}
where $(x,y)$ denote the usual  orthogonal Cartesian coordinates
parallel and normal to the boundary $ y = 0 $ (the wall),$ u(x,y),
v(x,y) $ are the $ x $ and
 $ y $ velocity components, respectively, and $ \nu > 0 $ is the kinematic viscosity.
   The function $ U_e=
U_e(x) $  is a given external velocity flow (the main-stream
velocity) which is assumed throughout the paper to be nonnegative
and is such that $ u(x,y) $ tends to $ U_e(x) $ as $ y\to \infty.$
 Equations (\ref{eq:Prandtl})  can be written in the form
\begin{equation}
\partial_y\psi\partial_{xy}^2\psi - \partial_{x}\psi\partial_{yy}^2\psi = U_e\partial_x U_e+\nu\partial^3_{yyy}\psi,
\end{equation}
where $ \psi $ is the well--known stream function defined by   $ u =
\ds\partial_{y}\psi, v = -\partial_{x}\psi.$

\noindent This equation with appropriate external velocity flow has
been the main focus of studies of particular exact solutions. It is
well known that to derive properties of solutions to a nonlinear
partial differential equations we use a family of special solutions.
They play an important role for describing the intermediate
asymptotic behavior of classes of solutions of original problems
with arbitrary initial data (see for example
\cite{Bar},\cite{Serrin}). A crucial step, in the analysis, is to
get favorable conditions such that particular solutions exist.  For
the boundary layer problems, research on this subject
 has a long history, which dates to the
pioneering works by  Blasius \cite{Blas} and Falkner and Skan
\cite{FS}. Their investigations lead to solutions to (1.2) in the
form
\begin{equation}\label{eq:similarity}
 \psi(x,y) = x^\alpha f(yx^{-\beta}).
\end{equation}
 Therefore, if
$ \psi $ is a such solution, it is easily verified that, for $ \beta
\not=0, $
\[ \psi(x,y)=\lambda^{-\frac{\alpha}{\beta}}\psi(\lambda^{\frac{1}{\beta}}x,\lambda y),\]
for all  $\lambda > 0, $ and we immediately see that
\[ \psi(x,y) = \left(\frac{x}{x_0}\right)^{\alpha}\psi\left(x_0,y\left(\frac{x}{x_0}\right)^{-\beta}\right).\]
This means that a solution $ \psi(x,y) $ for $ y $ fixed is similar
to the solution $ \psi(x_0,y) $ at a certain $ x_0.$
 This
solution is called invariant or similarity solution and the function
$ f $ is called the shape function or the dimensionless stream
function.\\  The main goal of identifying similarity solutions is to
reduce the original problem to an ordinary differential equation
which is easier to analyze. One says that the function (1.3) scales
the partial differential equation (1.2)  if the function $ f $
satisfies an ordinary differential equation called similarity
boundary layer equation. We refer the reader to
\cite{B7bis},\cite{BC},\cite{CMP6}, \cite{Ib1}, \cite{Ib2},
\cite{LCB} and the references therein. A relation of the type $
\alpha = h(\beta) $ called scaling relation.
\\ In \cite{Blas} the Blasius equation
\begin{equation}\label{eq:Blasius}
f^{\sec\prime}+ \frac{1}{2}ff^\sec = 0\quad \mbox{ on}\quad
(0,\infty)\end{equation} is derived from  (1.1) and (1.3)  where $
U_e $ is assumed to be  a constant function.  In the above equation
the primes indicate differentiations with respect to the similarity
variable $ t = yx^{-1/2}.$  Here we have $ \alpha = \beta =
\frac{1}{2}.$
 The Blasius equation is also a particular case  of the
Falkner-Skan equation \cite{FS}
\begin{equation}\label{eq:FS}
f^{\sec\prime}+ \frac{m+1}{2}ff^\sec = m({f^\prime}^2-1)\quad \mbox{
on}\quad (0,\infty),
\end{equation}
where  the external velocity  $ U_e(x) = U_\infty x^m, (U_\infty >
0).$ The exponents  $\alpha, \beta $  are given by $
\alpha=\frac{m+1}{2}, \beta
=\frac{1-m}{2}.  $  \\
Guided by the results of \cite{Blas}, \cite{FS} attention will be
given, in Section 2,   to identify a class  of  external velocities
for equation (1.2) to possess  solutions under the form (1.3) where
$ \alpha +\beta =1.$ Additionally, we investigate, in  Section 3,
the similarity solutions to (1.2), where both of the external
velocity and the stretching velocity of the permeable surface  are
assumed to vary as $ x^{-1}.$
\section{Similarity solutions}
\setcounter{equation}{0} \setcounter{theorem}{0}
\setcounter{lemma}{0} \setcounter{remark}{0}
\setcounter{corollary}{0} \qquad As it is said in the introduction,
this work deals with the similarity steady boundary--layer flow
induced by an incompressible viscous fluid  past a semi--infinite
flat plate. The phenomenon is governed by the system
\begin{equation}
u\partial_x u + v\partial_y u= \nu\partial^2_{yy} u +
U_e\partial_xU_{e}, \quad \partial_x u + \partial_y v = 0.
\end{equation}
accompanied by the  boundary condition
\begin{equation}
 u(x,y) \to U_e(x) \quad \mbox{ as } y \to \infty.
\end{equation}
 Therefore
 the stream function $\psi $ satisfies
\begin{equation}
\partial_y\psi\partial_{xy}^2\psi - \partial_{x}\psi\partial_{yy}^2\psi = U_e\partial_x U_e+\nu\partial^3_{yyy}\psi,
\end{equation}
and
\begin{equation}
 \lim_{y\to \infty}\partial_y\psi(x,y) = \partial_y\psi(x,\infty) = U_e(x).
\end{equation}

The main  problems, arising in the study of similarity solutions,
are related to the existence of the exponents $ \alpha $ and $ \beta
$ and to the rigorous study of the differential equation satisfied
by the profile $ f,$ which is, in general, nonlinear. For equation
(2.3), the classical approach for identifying $ \alpha $ and $ \beta
$ is the  scaling and transformation group \cite{Bar}. The essential
idea is to seek $ a $ and $ b $ such that if $ \psi $ satisfies
(2.3)   the new function $ \psi_\lambda(x,y) =
\lambda^a\psi(\lambda^b x, \lambda y) $ is also a solution.  This
will be certainly possible if the external velocity field $ U_e $ is
subject to the transformation. The parameters $ a $ and $ b $ may
depend on $ \alpha $ and $ \beta. $ Using (2.3), it is easily
verified that $ \psi_\lambda $  is a solution, for any $ \lambda, $
if the following
\begin{equation}\label{eq:conditionFS}
 a +b = 1,\quad \lambda^{3+a}U_e(\lambda^bx)U_{e_x}(\lambda^bx)=U_e(x)U_{e_x}(x),\end{equation}
hold for any $ x $ (see the proof of Theorem 2.1 below). In
\cite{FS} Falkner and Skan considered the case where
\begin{equation}\label{velocity} U_e(x) = U_\infty x^m,\end{equation}
where $ U_\infty(> 0) $ and $ m $ are constants. The case $ m = 0 $
was treated earlier by Blasius \cite{Blas}. Therefore, we deduce
from (2.5), for $ m \not=1,$
\[ a = \frac{m+1}{m-1},\quad b = -\frac{2}{m-1}.\]

Now  assume that  $ \psi $ satisfies the invariance property $
\psi_\lambda = \psi.$  Thus if we set $\lambda^bx = 1, $ we find
that $ \psi $ can be written in the form
\[ \psi(x,y) = x^{\alpha}f(y{x^{-\beta}}),\]
where
\[ \alpha = -\frac{a}{b}, \beta = \frac{1}{b}.\]
This leads to (1.3) with $ \alpha = \frac{m+1}{2} $ and $ \beta
=-\frac{m-1}{2}.$ This, may be, helps us to understand the
similarity stream functions obtained by Blasius and Falkner and
Skan. These authors considered solutions to (2.3),(2.4) of the form
\begin{equation}
\psi = (\nu U_e x)^{\frac{1}{2}}f(t),\quad t = (U_e/\nu
x)^{\frac{1}{2}}y.
\end{equation}
In fact, we shall see in the following results, that condition
(\ref{velocity}) is necessary for problem (2.3),(2.4) to admit
similarity solutions under the form (1.3). The main result of this
section says  that (1.3) scales  (2.3)
 if  the external velocity satisfies
\[ U_e(x) = U_\infty(x+x_0)^m,\]
for some $ x_0, $ where $ m = \alpha-\beta.$  We note, in passing
that if $ \psi $ is a solution the function $ \overline{\psi}(x,y) =
\psi(x+x_0, y)$ satisfies (2.3) with $ U_e(x+x_0) $ instead of $
U_e(x).$ Then it allows $ x_0 $ to be zero. The  necessary condition
on $ U_e $ can now be stated as follows.

\begin{Th}\label{Th:2.1} Assume that equation {\rm(2.3) } has a similarity solution in the form {\rm (1.3) } where $\alpha + \beta = 1.$ Then, there exist two
nonnegative constants, $c_1 , c_2 $  such that
\begin{equation}\label{eq:UE}
 U_e^2(x) = {c_1x^{2m}+c_2},\end{equation}
for all $ x > 0, $ where $  m = \alpha-\beta.$
\end{Th}
\begin{Proof}
 Let $ \psi $ be a stream function to (2.3) defined by (1.3) where $ \alpha + \beta = 1.$ Assume first that $ \beta\not=0.$  We choose $ a =
-\frac{\alpha}{\beta},\ b = \frac{1}{\beta},$ and define
$\psi_\lambda(x,y) = \lambda^a\psi\left(\lambda^bx,\lambda
y\right).$ Hence $ a + b = 1, \psi\equiv \psi_\lambda $ and
\[ L(\psi_\lambda)(x,y) = \lambda^{a+3}L(\psi)(\lambda^bx,\lambda y )\] for any $ \lambda > 0,$
where $ L $ is the operator defined by
\[ L(\psi) =\partial_y\psi\partial_{xy}^2\psi - \partial_{x}\psi\partial_{yy}^2\psi-\nu\partial^3_{yyy}\psi.\]
According to equation (2.3) we deduce for any $ \lambda > 0 $ and
any $ x > 0,$
\[ h(x)=\lambda^{a+3}h(\lambda^bx),\]
where $ h(x) =  U_e(x)\partial_xU(x).$ In particular, for fixed $
x_0
> 0 $
\[ h(\lambda^bx_0)=\lambda^{-(a+3)}h(x_0).\]
Setting  $ x = \lambda^bx_0 $ we infer
  \[h(x)=x^{-\frac{a+3}{b}}x_0^{\frac{a+3}{b}}h(x_0).\]
Solving the equation
\begin{equation} U_e\frac{d U_e}{dx} =  x^{-\frac{a+3}{b}}x_0^{\frac{a+3}{b}}h(x_0)\end{equation} yields us  (2.8) for $\beta\not=0$, since $ -\frac{a+3}{b}
+1=2(\alpha-\beta). $\\ For  $ \beta = 0, $ hence $ \alpha =m = 1,$
the new function
\[ \psi_\lambda(x,y)= \lambda^a\psi(\lambda^{-a}x,y),\]
for any fixed $ a \not= 0, $ is equivalent to $ \psi $ and satisfies
\[ L(\psi_\lambda)(x,y) = \lambda^{a}L(\psi)(\lambda^{-a}x,y )\] for any $ \lambda > 0.$ Arguing as in  the case $ \beta \not=0 $ one arrives at (2.8) with $ m
= 1.$
\end{Proof}

A similarity  assumption to (2.9) was proposed by Spalding \cite{Sp}
and studied extensively by Evans \cite{E}. The authors assumed that
$ U_e $ satisfies the equation
\begin{equation}\label{SE}\frac{dU_e}{dx} = C U_e^{2(\gamma-1)/\gamma},\end{equation}
where $ C $ is a constant. The similarity equation is then
\begin{equation}\label{eq:Cop} f^{\sec\prime} + ff^\sec + \gamma(1-{f^\prime}^2) = 0,\end{equation}
which is equivalent to (1.5) with $ \gamma = \frac{2m}{m+1}.$
  Now let us
discuss the consequence of the boundary condition at infinity. We
note, in passing, that from (\ref{eq:UE}) and (\ref{SE}) we easily
deduce that $ U_e(x) = U_\infty x^m.$
\begin{Th}
 Equation {\rm(2.3)} has a similarity solution,
$\psi,$  in the form {\rm(1.3)} such that $\alpha + \beta = 1$ and
$\partial_y\psi(.,\infty)=U_e,  $ if and only if there exists a real
$ m $ such that Then
\begin{equation}\label{eq:UEbis}
 U_e(x) = U_\infty x^{m},\end{equation}
for all $ x > 0, $  $ U_\infty $ is a constant. Moreover the real $
m $ satisfies  $  m = \alpha-\beta. $
\end{Th}
\begin{Proof} Let $ U_e(x) = U_\infty x^m.$ the existence of similarity solutions in the form (1.3) where $ \alpha = \frac{m+1}{2} $ and $
\beta=\frac{m-1}{2} $ is given in \cite{Blas} and \cite{FS} (see
also \cite{HT1},\cite{HT2}) for some values of $ m.$ Conversely,
assume that  (2.3), (2.4) has a similarity solution in the form
(1.3) with $ \alpha+\beta = 1.$  Then $ U_e(x)^2 =c_1x^{2m}+ c_2 $
and $ \lim_{y\to\infty}x^{2m}(f^\prime(yx^{-\beta})^2 = c_1x^{2m}+
c_2.$ Thereafter, the function $ {f^\prime}^2 $ has a finite limit
at infinity, which is unique and is given by $ c_1 + c_2x^{-2m}.$
This is acceptable only for $ c_2 = 0.$
\end{Proof}
\begin{Rm} {\rm It may be noted that in the above theorem it is not required to find the range of the real $ m $ such that problem (2.3),(2.4) has  a
similarity solution in the form (1.3). In fact, the similarity
equation may have no solution for some real $m $ (see Section 3).
Theorem 2.2 indicates, in particular, that for a prescribed external
velocity $ U_e(x) = U_\infty x^m,$ the reals $ \alpha $ and $ \beta
$ such that (1.3) scales (2.3) are given by }
\[\alpha =\frac{m+1}{2} , \beta =-\frac{m-1}{2} . \]
\end{Rm}
\begin{Rm}{\rm For a general external velocity, our approach can be used to obtain  particular solutions. In \cite{LCB} it is obtained particular solutions,
 having the form
\[ \psi = H(x,y)f(\eta(x,y)) + \psi_0.\]
So, if \[ U_e(x) = Ax^m + {\cal U}(x) \] we conjecture that there
exist solutions in the form
\[ \psi(x,y) = x^\alpha f(yx^{-\beta}) + \psi_0(x,y), \]
where $ \alpha = \frac{m+1}{2} $ and $ \beta = -\frac{m-1}{2}.$ The
function $ \psi_0 $ may be similarity and connected to $ x^\alpha
f(yx^{-\beta}) $ and $ {\cal U}(x). $
In the next result we shall identify the external velocity  such
that equation (2.3) has a solution in the form
\begin{equation}
 \psi(x,y)=x^\alpha f(yx^{-\beta}) + cx^ny,\ \alpha+\beta=1,\quad c = const.\end{equation}
Note that the function $ \psi_0(x,y) = cx^ny $ can be written in the
form $ \psi_0(x,y)= cx^{\alpha_1}(yx^{-\beta_1}), $ for any $
\alpha_1, \beta_1 $ such that $\alpha_1-\beta_1 = n; $ that is $
\psi_0 $  is similarity. So, for $ \alpha-\beta = n$ (2.13) is
identically to (1.3). Therefore we assume that $ \alpha
-\beta\not=n.$
 }
\end{Rm}

\begin{Th} Let $ \alpha, \beta $ are reals such that $ \alpha-\beta\not=n.$ The reals $ \alpha, \beta $ in {\rm(2.13) } scale {\rm(2.3)}
if and only if   $ \alpha =\frac{2}{3}, \beta =\frac{1}{3},\ n =
-\frac{1}{3} $ and the following
 \[ U_e(x) =Ax^{\frac{1}{3}} + cx^{-\frac{1}{3}},  \]
 holds for some constant $ A.$
\end{Th}
\begin{Proof} For $  U_e(x) = ax^{1/3} + cx^{-1/3},$ it is shown in \cite{LCB} that the function  $\psi(x,y) = x^{2/3}f(yx^{-1/3}) + cyx^{-1/3}$
satisfies (2.3) \\
Now, assume that (2.3) has a solution in the form (2.13).
Substituting  (2.13) into (2.3) yields
\begin{equation}\left\{ \begin{array}{rl}
x^{2\alpha-2\beta - 1}\Big[\nu f^{\sec\prime}\ +\ \alpha ff^\sec &
\left.-(\alpha-\beta){f^\prime}^2\right]\ +\
c(n+\alpha-\beta)x^{n+\alpha-\beta
-1}f^\prime\\
\\
 \ &+nc(n+\beta)tx^{n+\alpha-\beta -1}f^\sec -cnx^{2n - 1}+U_e\frac{d\ U_e(x)}{dx} = 0,\end{array} \right.\end{equation}
where $ t = yx^{-\beta}.$ Since $ \alpha-\beta\not = n, $ the
function $ f $ satisfies an ordinary differential equation if and
only if
\begin{equation}
 n + \alpha -\beta =n+\beta = 0, \end{equation}
and there exists a real $ \gamma $ such that
\begin{equation}\left\{ \begin{array}{lll}
\nu f^{\sec\prime}\ +\ \alpha ff^\sec  -(\alpha-\beta){f^\prime}^2&=&\gamma,\\
\\
  -cnx^{2n - 1}+U_e\partial_xU_e &=& \gamma x^{2(\alpha-\beta)-1}.\end{array} \right.\end{equation}
From (2.13) and (2.14)$_2$ we deduce  immediately  that $ \alpha
=\frac{2}{3}, \beta =\frac{1}{3},\ n = -\frac{1}{3} $ and $ U_e(x)
=Ax^{\frac{1}{3}} + cx^{-\frac{1}{3}}.$
\end{Proof}

\begin{Rm}{\rm Let us now  derive the well known   Blasius and Falkner-Skan.
Of course the external velocity is given by $ U_e(x) = U_\infty
x^m.$ We recall that  for the  Blasius model we have $ m = 0$ and
the case $ m\not= 0 $ was considered  by Falkner--Skan. Instead of
taking $ \alpha = \frac{m+1}{2}$ and $ \beta = -\frac{m-1}{2}$ we
shall insert (1.3) into (2.3) and choose $ \alpha $ and $ \beta $
such that $ f $ satisfies an ordinary differential equation.
Obviously we shall obtain that $ \alpha $  and $ \beta $ must to be , respectively  $\frac{m+1}{2}$ and $ -\frac{m-1}{2}. $ \\
Problem (2.3),(2.4) is written as
\begin{equation}
\partial_y\psi\partial_{xy}^2\psi - \partial_{x}\psi\partial_{yy}^2\psi =\nu\partial^3_{yyy}\psi+ mU_\infty^2x^{2m-1},
\end{equation}
with the boundary conditions
\begin{equation}
\partial_x\psi(x,0) = \partial_y\psi(x,0) = 0,\quad\partial_y\psi(x,\infty) = U_\infty x^m.
\end{equation}
If we substitute (1.3) into the first equation of this  problem we
obtain
\[ (\alpha-\beta)x^{2(\alpha-\beta)-1}(f^\prime)^2-\alpha x^{2(\alpha-\beta)-1}ff^\sec- x^{\alpha-3\beta}f^{\sec\prime}-mU_\infty^2x^{2m-1}=0,\]
and this is an ordinary differential equation if and only if
\[2(\alpha-\beta) - 1 = \alpha - 3\beta = 2m-1, \]
(the scaling relation) i.e.
\[ \alpha = \frac{m+1}{2}\quad\mbox{ and }\quad\beta = -\frac{m-1}{2}.\]
After a scaling ( $ \sqrt{\nu U_\infty}
f\left(.\sqrt{\frac{U_\infty}{\nu}}\right)$ instead of $ f(.) $) the
corresponding ordinary differential equation is
\begin{equation}
  f^{\sec\prime}+\frac{m+1}{2}ff^\sec + m(1-{f^\prime}^2)=0.\end{equation}
The prime denotes the derivative with respect to $ t =
\sqrt{\frac{U_\infty}{\nu}}yx^{\frac{m-1}{2}}. $ The boundary
conditions read
\begin{equation}
f(0) = f^\prime(0) = 0, f^\prime(\infty) = 1.
\end{equation}
}
\end{Rm}
\begin{Rm}{\rm We observe, in passing, that the  boundary conditions on the plate  are not required. This means that our analysis works even if we have a
continuous stretching surface. In this case,  if the stretching
velocity is given by $ U_s(x) = U_wx^n, $ we deduce from (1.3), $
n =\alpha-\beta $ which is the same exponent as that of the external
velocity. If  $ v = 0 $  at $ y = 0 $ the function $ f $ satisfies
\cite{HT1},\cite{HT2}
\[ f^{\sec\prime} + \frac{m+1}{2}ff^\sec -m({f^\prime}^2-1) = 0, \]
\[ f(0) = 0, \ f^\prime(0) = \zeta, \ f^\prime(0) = 1, \]
where $ \zeta $ is equal to the ratio of the free stream velocity to the boundary velocity.\\
In the case where the external velocity is zero, the stretching and
the suction/injection velocity have the form
\[ u(x,0)= U_wx^{m},\quad v(x,0) = -V_wx^{\frac{m-1}{2}},\]
the similarity equation is given by
\begin{equation}
\left\{\begin{array}{l}
 f^{\sec\prime} + \frac{m+1}{2}ff^\sec -m{f^\prime}^2 = 0,\\
\\
f(0) = a, \ f^\prime(0) =1, \ f^\prime(\infty) =
0,\end{array}\right.
\end{equation}
where $ a $ is  a real (the suction/injection parameter). This
problem arises also in the study the free convection, along a
vertical flat plate embedded in a porous medium.   We  refer the
reader to the papers
\cite{MK7},\cite{CMP6},\cite{BBT3},\cite{Guedda2} and the references
therein for a complete physical derivation and analysis of this
problem. }
\end{Rm}
\begin{Rm}{\rm We finally, mention a result of  \cite{Serrin} which provides  that for  the external velocity satisfying (2.12) where $
U_\infty > 0 $ and $ m > 0,$ any solution $ u $ to (1.1) with $
U_w=V_w = 0 $ such that $u_y $ is continuous for $ 0 < x < \infty, 0
\leq y < \infty, $ satisfies  the fundamental asymptotic estimate
\[ \left\vert \frac{u}{U_e}-f^\prime\right\vert = o\left(\frac{1+m\log x}{x^m}\right),\]
as $ x\to \infty,$ uniformly in $ y,$ where $ f $ is a solution to
(2.19),(2.20). }
\end{Rm}
\section{ The pseudo-similarity  solutions}
\setcounter{equation}{0} \setcounter{theorem}{0}
\setcounter{lemma}{0} \setcounter{remark}{0}
\setcounter{corollary}{0}\qquad In the present section we focuss our
attention to the case $m=-1$. The case $ m \not=-1 $ has been
abundantly studied. In particular it was considered for the  Falkner
Skan equation (2.11) where $ \gamma = 2m/(m+1).$ The cases $ \gamma
= 0 $ and $ \gamma = 1/2 $ are referred as Blasius and Homann
differential equations, respectively. In \cite{Cop11} Coppel
classified all solutions of (\ref{eq:Cop}), where $ 0 \leq \gamma <
2.$ Craven and Peletier \cite{CraPe} proved that equation (2.11),
where $\frac{1}{2} < \gamma \leq 1, $ has at most one solution $ f $
satisfying the boundary condition
\begin{equation}
 f(0) =f^\prime(0) = 0,\ f^\prime(\infty) = 1.\end{equation}
  In \cite{HS} Hasting and Siegel showed that for $\vert \gamma\vert
$  sufficiently small there exists a unique solution to (2.11),(3.1)
satisfying $ \vert f^\prime \vert < 1 $ and $ f^\sec(0) < 0.$ The
case $ \vert \gamma\vert > 1 $   can be found in the works
\cite{HT1,HT2} by Hasting and Troy. In \cite{HT1} deals with $\gamma
> 1.$ In this  it is shown that the equation has a periodic solution
and for any integer $ N $ problem (2.11),(3.1) has a solution with
at least $ N $ relative minima. It is also mentioned that as $
\gamma$ increases, the structure of periodic solutions and solutions
to (2.11),(3.1) gets progressively more complicated.
\\ In the present section we restrict our attention to the case $ m = -1.$ Therefore we shall consider the external velocity is given by $ U_e(x) =
\frac{U_\infty}{x} ( m =-1),$ and get new solutions to  the problem
 \begin{equation}
\partial_y\psi\partial_{xy}^2\psi - \partial_{x}\psi\partial_{yy}^2\psi = \nu\partial_{yyy}^3\psi-U_\infty x^{-3},
\end{equation}
subject to  the  boundary conditions
\begin{equation}
\partial_y\psi(x,0) = U_w.x^{-1},\quad \partial_x\psi(x,0) = -V_wx^{-1},\quad \partial_y\psi(x,\infty)=U_\infty.x^{-1},
\end{equation}
where  $ V_w $ is a real and  $ U_w $ and $ U_\infty $ are
nonnegative  and satisfy $ U_w <  U_\infty.$  The subject of the
present section is motivated by the work of Magyari, Pop and Keller
\cite{MPK},where the case $ U_\infty = 0 $ was considered and
discussed in detail. The authors  showed that if $  m = -1 $ problem
(2.21) has no solution. In order to overcome this difficulty  the
authors showed that the term $ V_{w}log(x) $ must be added to
expression (1.3). It has been also confirmed, by numerical
calculations, that  new solutions (pseudo-similarity solutions) exist provided that the suction parameter is large. \\
 For (3.2),(3.3), where $ U_\infty\not=0$ and according to Section 2, the function
$ \psi $ can be written as
\[ \psi(x,y) = \sqrt{\nu U_\infty}f(\sqrt{\frac{U_\infty}{\nu}}yx^{-1}).\]
Since, in general $ V_w = \frac{m+1}{2}\sqrt{\nu U_\infty}f(0),$ we
deduce  $ V_w = 0 $ for $ m=-1 $  and
  the function $\theta= f^\prime $ satisfies
\begin{equation}
\left\{\begin{array}{lll}
 \theta^{\sec}  +{\theta}^2-1 = 0, \\
\\
\theta(0) = \zeta, \theta(\infty) = 1,
\end{array}\right.\end{equation}
where $ \zeta = \frac{U_w}{U_\infty} $ is in the interval $ [0,1).$\\
The stability of equilibrium point  $ (1,0) $ of (3.4) cannot be
determined from the linearization. To analyze the behavior of the
nonlinear equation (3.4)$_1$, we observe that
\[ E^\prime(t)=0,\]
where $E$ is the Lyapunov function defined by
\[ E(t) = \frac{1}{2}\theta^\prime(t)^2 + \frac{1}{3}\theta(t)^3 - \theta(t).\]
Then, for  some  constant $ C,$ the following
\[ \theta^{\prime}=\pm\sqrt{2}\left(C+\theta-\frac{1}{3}\theta^3\right)^{1/2},\]
holds. The analysis of the algebraic equation of the phase path in
the phase plane $(\theta,\theta')$ reveals that the equilibrium
point $(1,0)$ is a center. Hence Problem (3.4) has no solution for
any $ \zeta >-1 $ except the trivial one $ \theta = 1,$ (see {\it
Fig.}$3.1$).

\begin{center}
\includegraphics[height=2.8in]{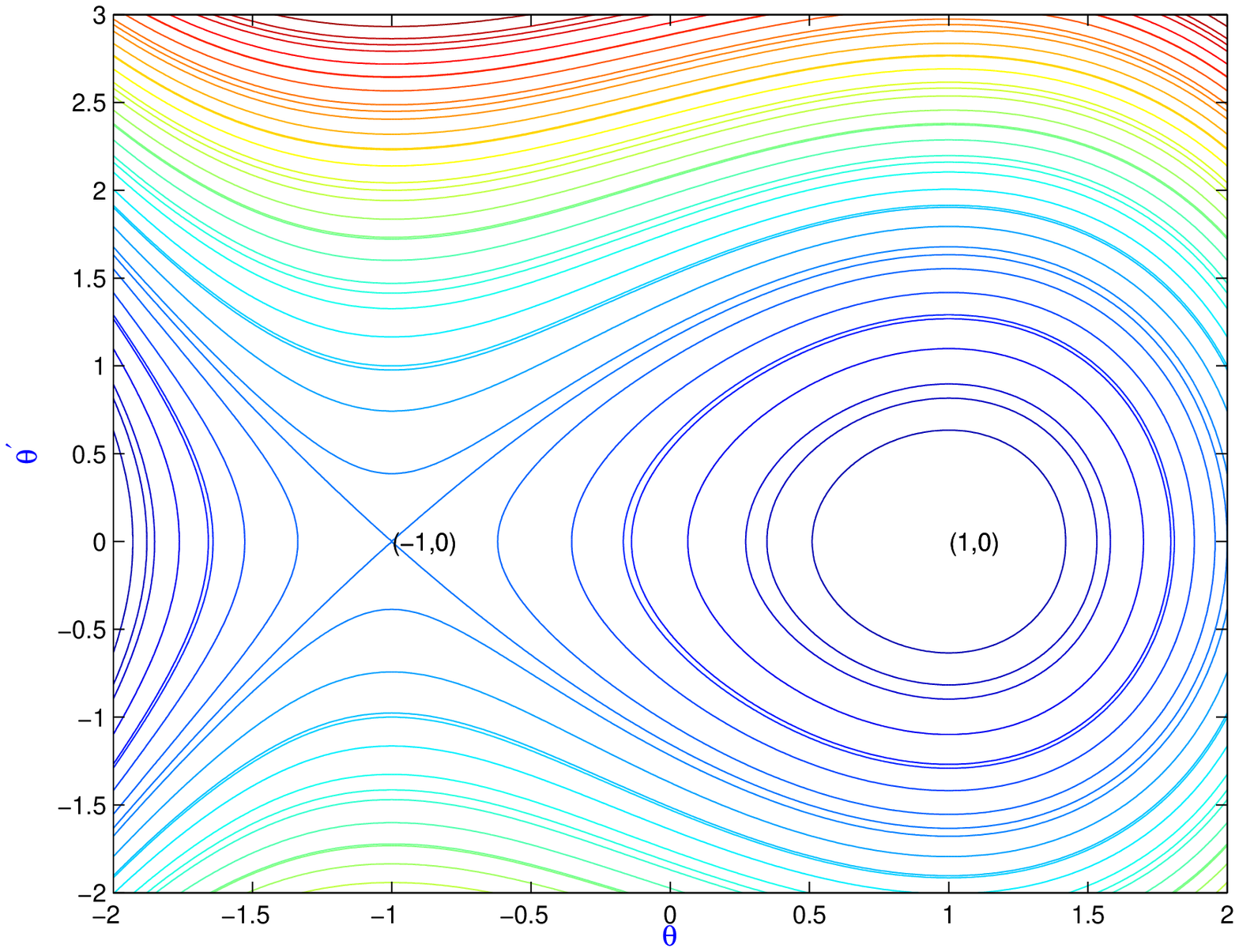}
\end{center}
\begin{center}
{\footnotesize
 {\it Fig. $3.1$} Classification of  solutions of $\theta^\sec + \theta^2 -1 = 0 $ according to $ \theta(0) $ and $ \theta^\prime(0).$}
\end{center}


\begin{Rm}{\rm If we impose the condition $ \theta(\infty) = -1 $ instead of $ \theta(\infty) = 1$--which is also of physical interest-- it is easy to see
that for any $ \zeta \leq 2 $ there exists a unique solution--up to
translation--. This solution satisfies
\[ \frac{1}{2}\theta^\prime(t)^2 + \frac{1}{3}\theta(t)^3 - \theta(t)= \frac{2}{3},\]
and we find that
\[ \theta(t) = 2-3\tanh^2\left[\pm t/\sqrt{2}+{\rm arctanh}\left\{(2-\zeta)/3)^{1/2}\right\} \right].\]
}
\end{Rm}
 \qquad To obtain exact solutions to (3.2),(3.3), we introduce the following Anzats for the stream-function
\begin{equation}\label{eq:Similarity}
 \psi(x,y) = ax^{\alpha} F(x,byx^{-\beta}).
\end{equation}
Guided by the analysis of Section 2 we take $\alpha
=\frac{m+1}{2},\beta = -\frac{m-1}{2},  a = \sqrt{\nu U_\infty} $
and $ b = \frac{U_\infty}{\nu}. $  Here the real $ m $ is
 assumed to be  any real. Hence boundary condition (3.3) reads
\[\partial_y\psi(x,0) = U_w x^{m},\quad \partial_x\psi(x,0) = -V_wx^{\frac{m-1}{2}},\quad \partial_y\psi(x,\infty)=U_\infty x^{m}.\]
Using this and the equation of the stream function
 we find that $ F $ satisfies
\begin{equation}
\left\{\begin{array}{l} F^{\sec\prime} +  \frac{1+m}{2} FF^{\sec}-
m({F^\prime}^2-1)+x\left(F^\sec\partial_xF-F^\prime\partial_xF^\prime\right)=0,\\
\\
F^\prime(x,0) = \zeta,\quad F^\prime(x,\infty)=1,\\
\\
\ds \frac{1+m}{2}F(x,0) +x\partial F(x,0) = \frac{V_w}{\sqrt{\nu
U_\infty}},
\end{array}\right.
\end{equation}
where the primes denote partial differential with respect to $ t =
\sqrt{\frac{U_\infty}{\nu}}yx^{\frac{m-1}{2}}.$ By writing
\[ F(x,t) = f(t) + H(x),\]
we find
\begin{equation}
\left\{\begin{array}{l}
 f^{\sec\prime} +  \frac{1+m}{2} ff^{\sec}- m\left({f^\prime}^2-1\right)+f^\sec\left(xH^\prime + \frac{1+m}{2}H\right)=0,\\
\\
\ds \frac{1+m}{2}f(0) +\frac{1+m}{2}H(x) + xH^\prime(x)=
\frac{V_w}{\sqrt{\nu U_\infty}},\quad f^\prime(0) = \zeta \in
[0,1),\quad f^{\prime}(\infty) = 1.
\end{array}\right.
\end{equation}
Hence, there exists a real $ \t $ such that
\begin{equation}
\left\{\begin{array}{l}
f^{\sec\prime} +  \frac{1+m}{2} ff^{\sec}- m\left({f^\prime}^2-1\right)+\t f^\sec=0,\quad t > 0,\\
\\
xH^\prime + \frac{1+m}{2}H = \t, \quad x > 0,\\
\\
\ds \frac{1+m}{2}f(0) +\t= \frac{V_w}{\sqrt{\nu U_\infty}},\quad
f^\prime(0) = \zeta,\quad f^{\prime}(\infty) = 1.
\end{array}\right.
\end{equation}
 First, let us note that (2.19) and (3.8)$_1$ are equivalent for $ m \not=-1.$ In this case the general solution of (3.8)$_2$ is
\[H(x) = Cx^{-\frac{1+m}{2}}+ \frac{2\t}{1+m}, \]
where $ C $ is a constant. Thus the stream function is given by  $
\psi(x,y) = aC + ax^{\frac{1+m}{2}}\left(f(t) +
\frac{2\t}{1+m}\right),$
and the new function $ g = f + \frac{2\t}{1+m} $ satisfies equation (2.19).\\
Thereafter, we will assume that $ m = -1$ and this leads to
\begin{equation}
\left\{\begin{array}{l}
f^{\sec\prime}+ \t f^\sec+ {f^\prime}^2-1=0,\\
\\
f^\prime(0) = \zeta,\quad f^{\prime}(\infty) = 1,
\end{array}\right.
\end{equation}
and
\[ H(x) = \t\log x + C,\quad C = const.\]
Then the required exact solution has the form
\begin{equation}
\psi(x,y) =  af(byx^{-1}) + a\t\log(x).
\end{equation}
This formula is similar to the one studied in \cite{KV}  in the
context of rough surface growth and can be regarded as solution with
dynamic scaling \cite{KV}. In passing, we note that from (3.8) one
sees $ \t = V_w\left(\nu U_\infty\right)^{-1/2}.$ This means that
the constant $ \t $ plays the role of
suction/injection parameter.\\
To study (3.9) it is  more convenience to analysis the second
ordinary differential equation
\begin{equation}
\left\{\begin{array}{l}
\theta^\sec+ \t \theta^\prime+ \theta^2-1=0,\\
\\
\theta(0) = \zeta,\quad \theta(\infty) = 1,
\end{array}\right.
\end{equation}
where $0 \leq \zeta < 1 $ and $ \t\not=0.$ In fact  the real $ \t $
will be taken in $ (0,\infty).$ The existence of solutions to (3.11)
will be proved by means of shooting  method. Hence, the boundary
condition at infinity is replaced by the condition $
\theta^\prime(0) = d,$ where $ d $ is a real. For any $ d $ the new
initial--value problem has a unique local solution $ \theta_d $
defined in its maximal interval of existence $(0,T_d), T_d \leq
\infty.$ We shall see that for an appropriate $ d $ the solution
 $ \theta_d $ is global and satisfies
\begin{equation} \theta_d(\infty) = 1.\end{equation}
 A simple analysis in the phase plane  reveals that  problem (3.11)   may have solutions only for  $\t > 0. $ In fact the ordinary differential equation in
(3.11) is considered as a nonlinear autonomous system in $
\mathbb{R}^2, $ with the unknown $(\theta,\theta^\prime),$ mainly
\begin{equation}
\left\{\begin{array}{l}
\theta^\prime=\varphi,\\
\\
\varphi^\prime =-\t \varphi + 1 - \theta^2,\\
\end{array}\right.
\end{equation}
subject to the boundary condition
\begin{equation}
\theta(0) = \zeta,\quad \varphi(0) = d.
\end{equation}
The linear part of the above system at the equilibrium point $(1,0)$
is  defined by the matrix
\[ J=\left(\begin{array}{ccc}
                                    0 & 1\\
\\
-2 & -\t
\end{array}\right). \]
The eigenvalues of $ J  $  are
\[ \lambda_1 = \frac{-\t -\sqrt{\t^2-8}}{2},\quad \lambda_2 = \frac{-\t +\sqrt{\t^2-8}}{2},\]
if $ \t \geq \sqrt{8} $ and for $  0 < \t  < \sqrt{8}, $
\[ \lambda_1 = \frac{-\t -i\sqrt{8-\t^2}}{2},\quad \lambda_2 = \frac{-\t +i\sqrt{8-\t^2}}{2}.\]
Therefore, the hyperbolic equilibrium point $(1,0) $ is
asymptotically stable if $ \t $ is positive and unstable for $ \t $
is negative. In particular problem (3.11) has no nontrivial
solutions if $ \t < 0.$ If $ \t > 0 $ we deduce from the above that
there exists $ \delta > 0 $ such that for any $ d $ and $ \zeta $
such that $  d^2  + (\zeta - 1)^2 < \delta^2 $ the local solution $
\theta_d $ is global and satisfies (3.12). In the following we
construct solutions to (3.11) where the condition  $  d^2  + (\zeta
- 1)^2 < \delta^2 $ is   not necessarily required. For a
mathematical consideration the parameter $ \zeta $ will be taken in
$(-1,\sqrt{3}].$  The following theorem deals with  nonnegative
values of $ \zeta.$
\begin{Th} Let $ 0 \leq \zeta \leq \sqrt{3} $ and $ d $ be a real such that
\begin{equation}
d^2\leq 2\zeta(1-\frac{\zeta^2}{3}).
\end{equation}
Then the local solution $ \theta_d $ is global, nonnegative and
tends to $ 1 $ as $ t $ approaches infinity.
\end{Th}
The key  point of the proof of this theorem is to find $ d $ such
that $ \theta_d(t) > 0 $ for all $ t > 0.$ To show this we consider
again the function $ E(\theta(t),\varphi(t))  =
\frac{1}{2}\varphi(t)^2 + \frac{1}{3}\theta^3-\theta.$ Along an
orbit we have
\[ \frac{d}{dt}E(\theta(t),\varphi(t)) = -\t\varphi(t)^2\leq 0.\]
Hence
\[  E(\theta_d(t),\theta_d^\prime(t)) \leq E(\zeta,d),\]
as long as $ \theta_d(t) $ exists; that is $ t < T_d.$  The
following result shows that $ \theta_d \geq 0 $ on $(0,T_d)$ and
then $ T_d = \infty.$
\begin{Lm} Let $ 0\leq \zeta \leq \sqrt{3} $ and $ d $ satisfying (3.15). Then $ \theta_d $  is nonnegative,
global, bounded and its first derivative also is bounded.
\end{Lm}
\begin{Proof} First we note that  from (3.13) and (3.14), there exists $ t_0 > 0,$ small, such that $ \theta_d $ is positive on $(0,t_0).$  Assume that $
\theta_d $ vanishes at some $ t_1 > t_0 $  and suppose that  $
\theta_d^\prime(t_1) \not= 0. $ Because
\[E(\zeta,d) \geq E(\theta_d(t),\theta_d^\prime(t)) \geq \frac{1}{2}\theta_d^\prime(t_1)^2,\]
for all $ 0 \leq t\leq t_1. $ we deduce $ \frac{1}{2}d^2 >
\zeta(1-\frac{1}{3}\zeta^2), $ which contradicts (3.15). Therefore $
\theta^\prime_d(t_1)= 0.$ In this case we deduce from the equation
of $ \theta_d $ that $ \theta_d^\sec(t_1) = 1 $ and then $ \theta_d
$  is nonnegative on a some neighborhood of $ t_1.$ Consequently the
local solution is nonnegative as long as there exists.  To show that
$ \theta_d $ is global we note that
\[ E(\zeta,d) \geq \frac{1}{2}\theta_d^\prime(t)^2+\frac{1}{3}\theta_d^3(t)-\theta_d(t) \geq -\frac{2}{3},\]
for all $ t \leq T_d$, since $ \theta_d $ is nonnegative. Hence $
\theta_d $ and (then) $ \theta_d^\prime   $ are bounded.
Consequently $ \theta_d $ is global. The lemma is proved
\end{Proof}

Now, we are ready  to prove Theorem 3.1. More  precisely we have.
\begin{Prop} Let $ \t > 0, \zeta \in [0,\sqrt{3}] $ and $ d $ such that {\rm (3.15)} holds.
Then, the global solution $ \theta_d $  to {\rm (3.13),(3.14)}
satisfies the boundary condition {\rm (3.12)} at infinity.
\end{Prop}
\begin{Pf}
On account of Lemma 3.1 we need only to show that $ \theta_d $ goes
to 1 as $ t $ approaches infinity.
Since this result is the broad goal of the present section, we give two proofs.  \\
{ \it The first proof.}

Using the fact that $ \theta_d $ and $ \theta_d^\prime $ are bounded
we deduce  from the equation of $\theta_d $ that $ \theta_d^\sec $
is bounded. On the other hand, since $ \frac{d}{dt}
E(\theta_d,\theta_d^\prime) = -\t{\theta_d^\prime}^2,$ the function
$ \theta_d^\prime $ is square integrable. Now, we use the identity
\[\theta_d^\prime(t)^3 = \theta_d^\prime(0)^3 + 3\int_0^t\theta_d^\prime(s)^2\theta_d^\sec(s)ds,\]
to show that $ \theta_d^\prime(t) $ has a finite limit as $ t $
tends to infinity  and this limit is zero. Next, we get, by
differentiation (3.11)$_1$
\begin{equation}
 \theta_d^{\sec\prime} = -\t \theta^\sec -2\theta_d\theta_d^\prime.\end{equation}
Multiplying equation (3.16) by $ \theta_d^{\prime\prime},$
integrating the equation obtained over $(0,t),$ we get
\[ \t \int_0^t\theta_d^\sec(s)^2ds = \frac{1}{2}\theta_d^\sec(0)^2+\theta_d(0)\theta_d^\prime(0)^2- \frac{1}{2}\theta_d^\sec(t)^2 -
\theta_d(t)\theta_d^\prime(t)^2 + \int_0^t\theta_d^\prime(s)^3ds.\]
This implies $\theta_d^\sec\ \in L^2(0,\infty).$ Since  $
\theta_d^{\sec\prime} $ is bounded, by using (3.16), one sees $
\lim_{t\to\infty}\theta_d^\sec(t) = 0.$ Finally, we deduce from
(3.11)$_1$ that $ \theta_d(t) $  goes to 1 as $ t
$ tends to infinity since $ \theta_d $ is nonnegative. \\
{ \it The second proof.}

This proof uses the Bendixson Criterion. Let  ${\cal T}$ be the
trajectory of $(\theta_d,\theta_d^\prime)$  in  the phase plane $
(0,\infty)\times\mathbb{R} $ for $ t \geq 0 $ and let
$\Gamma^+({\cal T})$ be its $w$-limit set at $+\infty.$  From the
boudedness of ${\cal T} $  it follows that $ \Gamma^+({\cal T})$ is
a nonempty connected and compact subset of
$(0,\infty)\times\mathbb{R}$ ( see, for example \cite[p
226]{Amann}). Moreover $ (-1,0) \not\in \Gamma^+({\cal T}),$ since $
\theta_d $ is nonnegative.  Note that if $\Gamma^+({\cal T})$
contains the equilibrium point $(1,0) $ then $\Gamma^+({\cal T}) =
\left\{(1,0)\right\},$ since $(1,0)$ is asymptotically stable.
Assume  that $(1,0)\ \not\in\Gamma^+({\cal T}),$
 Applying  the Poincar\'e--Bendixson Theorem \cite[p 44]{Holmes} we deduce that $\Gamma^+({\cal T})$ is a
cycle, surrounding $(1,0)$. To finish the second proof we shall
prove the nonexistence of such a cycle. To this end we define $
P(\theta,\varphi) = \varphi, Q(\theta,\varphi) = -\t\varphi
+1-\theta^2, \varphi=\theta_d^\prime $ and $ \theta = \theta_d.$
 The function $(\theta,\varphi)$ satisfies the system $ \theta^\prime = P(\theta,\varphi),\ \varphi^\prime = Q(\theta,\varphi).$ Let $ D $ be the bounded
domain of the $(\theta,\varphi)$--plane with boundary ${\Gamma^+}. $
As $ P $ and $ Q $ are  regular  we deduce, via the Green--Riemann
theorem,
\begin{equation}
\int\int_D \left(\partial_\varphi Q + \partial_\theta
P\right)d\varphi d\theta=\int_{{\Gamma^+}}\left(Qd\theta - P
d\varphi\right)=0,
\end{equation}
thanks to  the system satisfied by $(\theta,\varphi). $  But $
\partial_\varphi Q + \partial_\theta P = \t $ which is positive. We
get   a contradiction.
\end{Pf}

In the same way as in the proof of Proposition 3.1, we can see that
any global solution to (3.13) which is bounded from bellow by some $
a > -1 $ tends to 1 at infinity. Therefore, to complete our
analysis, we shall determine a domain  of attraction of the critical
point $(1,0).$   Let
\[{\cal P} = \left\{(\zeta,d)\ \in \mathbb{R}^2: \zeta > -1, \frac{1}{2}d^2 + \zeta\left(\frac{1}{3}\zeta^2-1\right) < \frac{2}{3} \right\}.\]
\begin{Prop} For any $(d,\zeta) $ in ${\cal P}$ the local solution to (3.13),(3.14) is global and converges to 1 at infinity.
\end{Prop}
\begin{Proof} Let us consider a one-parameter
of family of curves defined by
\[ E(\theta,\varphi) =  \frac{1}{2}\varphi^2 + \frac{1}{3}\theta^3-\theta = C,\]
where $ C $ is a real parameter. Note that, in the phase plane, this
family is  solution curves of system (3.4).  The curve $ \varphi^2 =
2\theta - \frac{2}{3}\theta^3 + \frac{4}{3},$ corresponding to $ C =
\frac{2}{3},$ goes through the point $(2,0)$ and has the saddle
$(-1,0) $ ($\t=0$) as its $ \alpha$ and $w$-limit sets, see Fig 3.2.
We denote  this solution curve by $ {\cal H},$ which is, in fact, an
homoclinic orbit and define a separatrix cycle for (3.4). We shall
see that the bounded open domain with the boundary $ {\cal H}$ is an
attractor set for $(1,0)$ of system (3.13) where $
 \t > 0.$  This domain is given by $ E(\theta,\varphi) = C, \theta > -1,  $ for all $ -\frac{2}{3} \leq C < \frac{2}{3}, $ which is $ {\cal P}.$ As $
\frac{d}{d t} E \leq 0 $ any solution, with initial data in $ {\cal
P} $ cannot leave $  {\cal P},$ (see the proof of Theorem 3.1 ). By
LaSalle invariance principle we deduce that   for any $(\zeta,d) $
in $  {\cal P } $ the $w$-limit set, $ \Gamma^+{(\zeta,d)} $ is a
nonempty, connected subset of $  {\cal P }\cap \left\{\varphi =
0\right\},$ (see \cite[p 234]{Amann}.) However, if $ \theta\not=1,
\varphi = 0 $ is a transversal of the phase--flow, so the $w$--limit
set is  $ \left\{(1,0)\right\}.$ This means that  $ {\cal P} $ is a
domain of attraction of the critical point $(1,0). $
\end{Proof}
\begin{center}
\includegraphics[width=6.5cm]{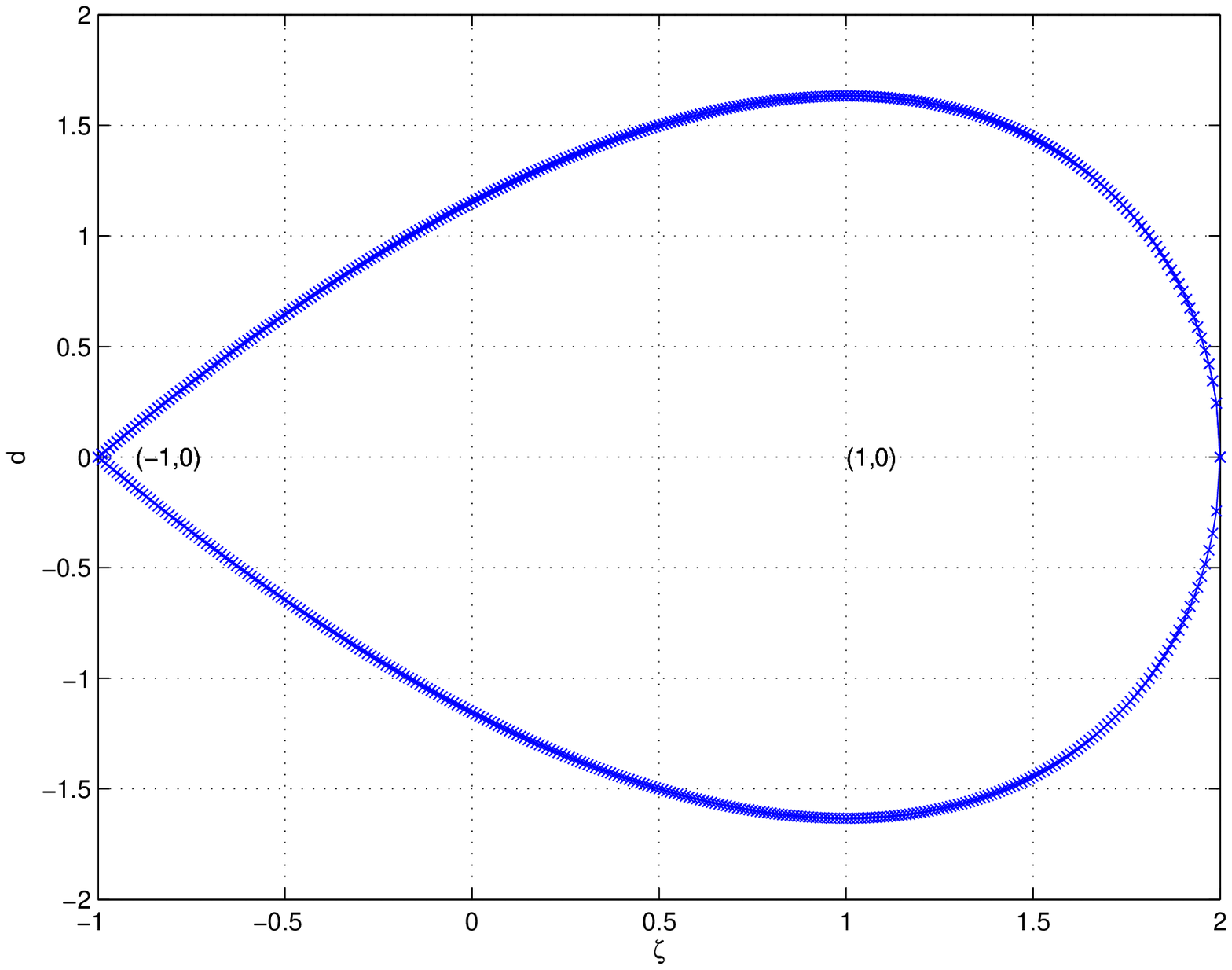}
\end{center}
\begin{center}
{\footnotesize
 {\it Fig. $3.2$} A basin of attraction ${\cal P }$ of the critical
point $(1,0)$.}
\end{center}
\section{Numerical results}
\qquad In this section numerical solutions of the boundary--value
problem $(3.11)$ are obtained by using the fourth-order Runge-Kutta
integration scheme with the shooting method. velocity profiles of
the dimensionless velocity $\theta$ are plotted in term of the
similarity variable $t$, for various value of the shooting parameter
$d$ (the dimensionless skin-friction), {\it Fig.}4.1 $(1)$ and $(2)$
show that the numerical results are in good agreement with the above
theoretical predictions.
\begin{center}

\includegraphics[width=6.5cm]{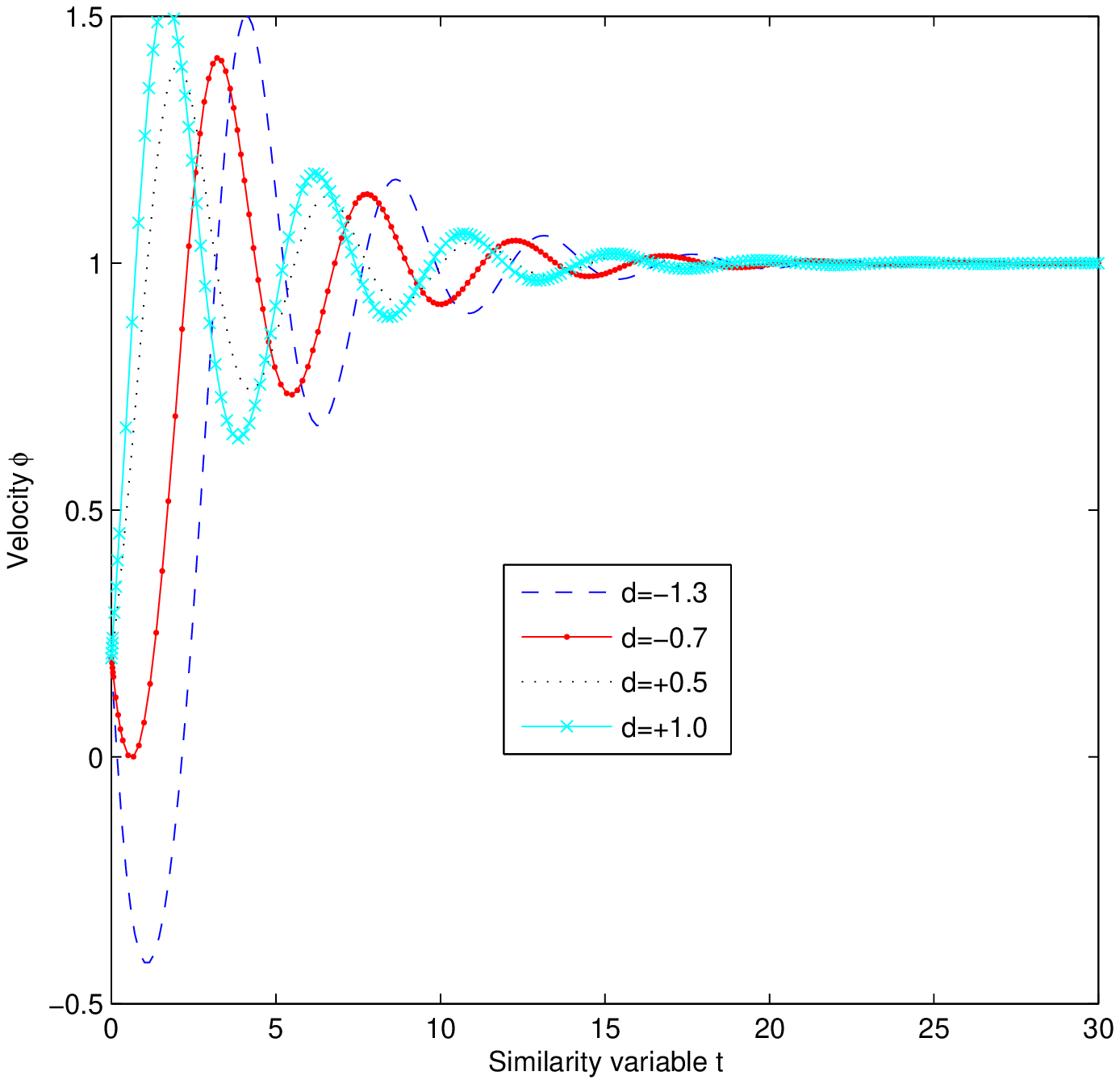}(1)
 
\end{center}
\begin{center}
\includegraphics[width=6.5cm]{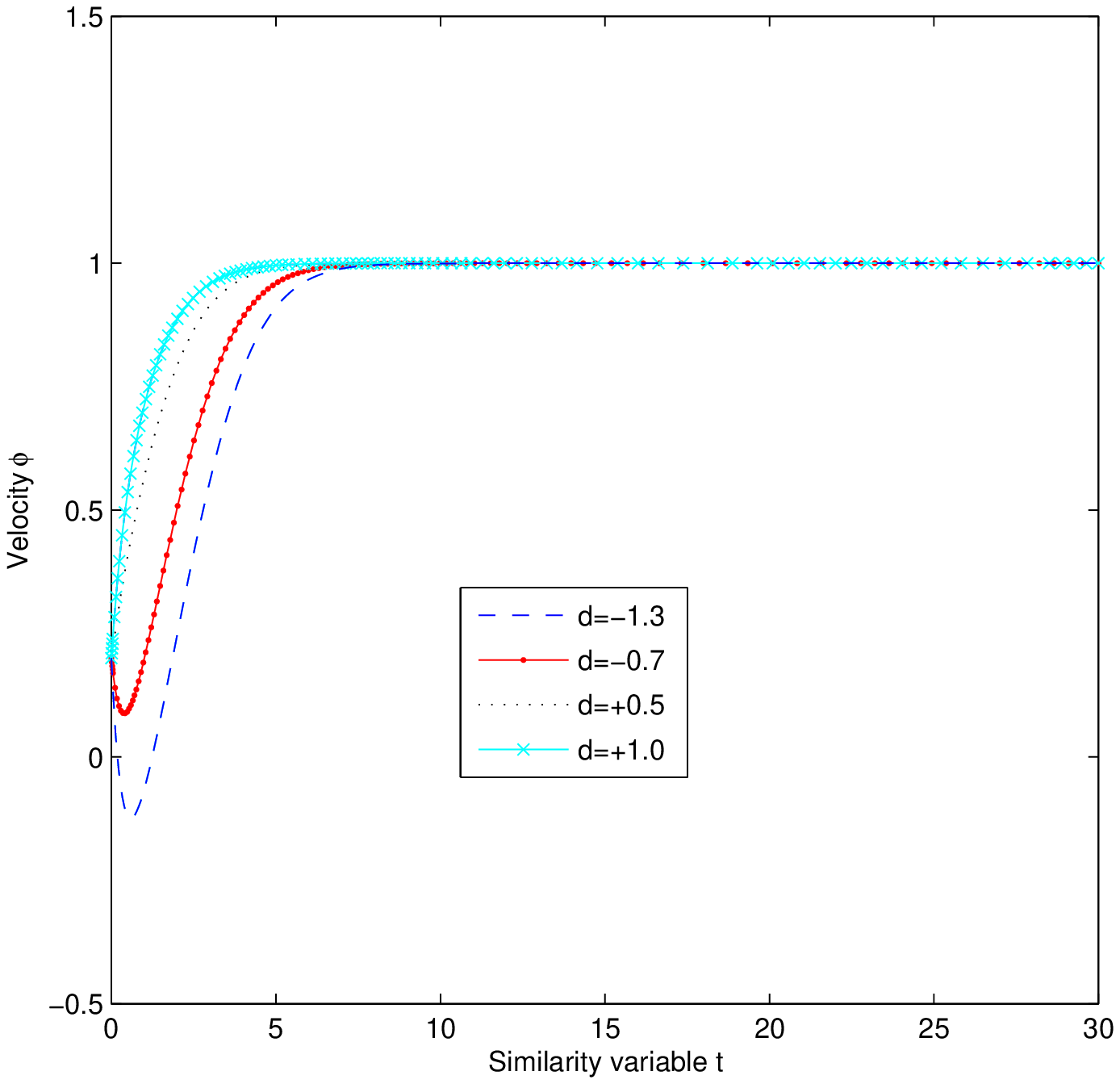}(2)

 \end{center}
\begin{center}
{\footnotesize {\it Fig. $1.4$} Velocity profiles in terms of
$d=\theta^{\prime}(0)$ for fixed $\zeta=0.2$ and $\tau=0.5$ $(1)$,
and $\tau=2.83$ $(2)$.}
\end{center}


\section{Conclusion}
\setcounter{equation}{0} \setcounter{theorem}{0}
\setcounter{lemma}{0} \setcounter{remark}{0}
\setcounter{corollary}{0} \qquad In this work the laminar
two-dimensional steady incompressible, boundary layer flow past a
stretching surface is considered.  It has been shown that the
problem has solutions having a similarity form if the velocity
distribution outside the boundary layer is proportional to $ x^m,$
for some real number $ m.$ In the second part of this paper, we are
interested in question of existence of solutions t in the case where
the external velocity is an inverse-linear function; $m=-1.$ This
situation occurs in the case of sink flow.  To obtain exact
solutions the stream function $ \psi $ is written under the form
\begin{equation}
\psi(x,y) = \sqrt{\mu U_\infty}f(t) + V_w\log(x).
\end{equation}
  It is shown that
the
  ordinary differential equation satisfied by $ f $  has multiple solutions for any $ V_w $ positive (suction) and no solution
  can exist if $ V_w
\leq 0$ (injection). A sufficient condition for the existence is
derived:
\begin{equation}\zeta > -1,\quad \quad \frac{1}{2}f^{\sec}(0)^2 + \zeta\left(\frac{\zeta^2}{3}-1\right) < \frac{2}{3}
.\end{equation}  We have obtained two family of solutions according
to  $ \tau = V_w\left(\mu U_\infty\right)^{-1/2}. $  If $ \t \geq
\sqrt{8} $, $f^\prime $ is monotonic  and goes to unity at infinity,
but if $ 0 < \t < \sqrt{8},$ we have a stable spiral.  The function
$f^\prime $ oscillates an infinite number of times and goes to 1. So
if we are interested in solutions to  $(3.9)$ such that
\[ -1 < f^\prime < 1 \] we must take $ U_w, V_w$ and $  U_\infty > 0 $
 satisfying $ -U_\infty < U_w < U_\infty $  and $ V_w > (8\nu
U_\infty)^{1/2}.$\\
Condition (5.2) indicates also that for the same positive value of
the suction parameter the permeable wall stretching with velocity $
U_wx^{-1}, U_w > 0 $ has multiple boundary--layer flows. Every flow
is uniquely determined by the dimensionless skin friction $
f^\sec(0) $ which can be any real number  in the interval
\begin{center}
$\left(-\sqrt{\frac{4}{3} +
2\zeta\left(1-\zeta^2/3\right)},\sqrt{\frac{4}{3} +
2\zeta\left(1-\zeta^2/3\right)}\right),$ \end{center} where $ \zeta
= U_wU_\infty^{-1}.$ The case $ U_\infty = 0 $ was considered by
Magyari, Pop and Keller \cite{MPK}. The authors showed, by numerical
solutions, that the boundary
layer flow exists only  for a large suction parameter ($\tau \geq 1.079131$).\\
 The existence of exact solutions of the Falkner-Skan equation under the
present condition was discussed by Rosenhad \cite[pp.
244--246]{Rose} who mentioned that these results may be obtained by
rigorous arguments which, in fact, motivated the present work.   We
note, in passing, that it is possible to obtain solutions if the the
skin friction satisfies $ \vert f^\sec(0)\vert > \sqrt{\frac{4}{3} +
2\zeta\left(1-\zeta^2/3\right)}.$

\end{document}